\documentclass[]{amsart}
 \usepackage[]{amsmath, amsthm, amsfonts,amssymb}
 \usepackage[all]{xy}
 \usepackage{hyperref}
 \xyoption{curve}
 \usepackage[mathscr]{eucal}

 \newcommand {\C} {{\mathbb C}}
 
 \newcommand {\Z} {{\mathbb Z}}
 \newcommand {\Q} {{\mathbb Q}}
 \newcommand {\PP} {{\mathbb P}}
 
 \newcommand {\F} {{\mathcal F}}

 \newcommand {\dt} {{\bullet}}

 \newcommand {\codim} {\text{codim}}
 \newcommand {\im} {\text{im}}

 \newcommand {\tX} {\tilde{X}}
 \newcommand{\V}{\mathscr{V}ar}
 \newcommand{\HS}{\mathscr{HS}}
 \newcommand{\FHS}{\mathscr{FHS}}
 
 \newcommand{\K}{\text{K}_0} 
 \newcommand{\mO}{\mathcal{O}}
\newcommand{\LL}{\mathbb{L}}
 \newcommand{\KK}{\mathbb{K}}
 \newcommand{\A}{\mathbb{A}}
 \newcommand{\M}{\mathcal{M}}
 \newcommand{\hM}{\widehat{\mathcal{M}}}
 
 \newcommand{\hN}{\widehat{\mathcal{N}}}
 
\newcommand{\cL} {\mathcal{L}}
\newcommand{\cC} {\mathscr{C}}
\newcommand{\FcC} {\mathscr{FC}}
\newcommand{\GcC} {\mathscr{GC}}

 \newtheorem{thm}[subsection]{Theorem}
 \newtheorem{cor}[subsection]{Corollary}
 \newtheorem{lemma}[subsection]{Lemma}

 \newtheorem{rmk}[subsection]{Remark}

\begin{document}

 \title{Coniveau and the Grothendieck group of varieties}

 \author{Donu Arapura}\thanks{First author partially supported by the
   NSF} \author{ Su-Jeong Kang}
 \address{Department of Mathematics\\
   Purdue University\\
   West Lafayette, IN 47907\\
   U.S.A.}  \maketitle

 There are two natural filtrations on the singular cohomology of a
 complex smooth projective variety: the coniveau filtration which is
 defined geometrically, and the level filtration which is defined
 Hodge theoretically. We will say that the generalized Hodge
 conjecture (GHC) holds for a variety $X$ if these filtrations
 coincide on its cohomology. There are a number of intermediate forms
 of this condition, including the statement that the ordinary Hodge
 conjecture holds for $X$. We show that if GHC (or an intermediate
 version of it) holds for $X$ then it holds for any variety $Y$ which
 defines the same class in a completion of Grothendieck group of
 varieties. In particular, using motivic integration we see that this
 is the case if $X$ and $Y$ are birationally equivalent Calabi-Yau
 varieties, or more generally $K$-equivalent varieties.  This refines
 a result obtained in \cite{Arapura} by a different method.

 The key point is to show that the singular cohomology with its
 coniveau (respectively level) filtration determines a homomorphism
 $\nu$ (respectively $\lambda$) from the Grothendieck group of
 varieties $\K(\V_{\C})$ to the Grothendieck group of polarizable filtered Hodge
 structures $\K (\FHS)$. This is done by showing that cohomology
 together with these filtrations behave appropriately under blow ups.
 We then show that $X$ satisfies GHC if and only if its class $[X]$
 lies in the kernel of the difference $\nu - \lambda$, and the above
 results follow from this.

 The following conventions will be used throughout the paper. All our
 varieties will be defined over $\C$. We denote the singular
 cohomology of a smooth projective variety $X$ with rational
 coefficients by $H^i(X)$. Our thanks to the referee for a number of
 helpful suggestions.

 \section{Filtered Hodge structures}

 Let $X$ be a smooth projective variety. Its cohomology carries a
 natural Hodge structure. The {\em coniveau} filtration on $H^i(X)$ is
 given by $$
 N^p H^i(X) = \sum_{\codim S\ge p} \text{ker}[H^i(X)\to
 H^i(X-S)]$$
 It is a descending filtration by sub Hodge structures.
 The largest rational sub Hodge structure $\F^pH^i(X)$ contained in
 $F^pH^i(X)$ gives a second filtration, which we call the level
 filtration.  We have $N^pH^i(X)\subseteq \F^pH^i(X)$. We will say
 that $\text{GHC}(H^i(X),p)$ holds if $N^pH^i(X)=\F^pH^i(X)$.  We will
 say that generalized Hodge conjecture (GHC) holds for $X$ if we have
 equality for all $i$ and $p$.  Note that $\text{GHC}(H^{2p}(X),p)$ is
 just the usual Hodge conjecture.

We recall that a Hodge structure is polarizable if it
admits a polarization, that is a bilinear form satisfying the
Hodge-Riemann bilinear relations. We note  that the Hodge structure on the
cohomology of a smooth projective variety is  polarizable:
once an ample line bundle is chosen, 
a polarization is  given by taking the orthogonal direct sum
of the polarizations, determined in the usual way,
 on  primitive cohomology \cite[p. 202, 207]{wells}.
 Let $\HS$ be the category of finite direct sums of pure rational
 polarizable Hodge structures. 
The category $\HS$ is a semisimple  Abelian category with tensor
products \cite[4.2.3]{deligne}.
 Any object $H$ can be decomposed into a sum $$H= \bigoplus_i H^i$$
 where $H^i$ is the largest sub Hodge structure of weight $i$. We
 define the category $\FHS$ whose objects are polarizable Hodge structures with
 finite descending filtrations by sub Hodge structures, and whose
 morphisms preserve the filtration. Note that this would be  a filtration
by polarizable sub Hodge structures, since  a sub Hodge structure of
 a polarizable Hodge structure is again polarizable.

Given an additive category $\cC$, we can define a Grothendieck
group $\K^{split}(\cC)$  by generators and relations as follows.   
We have one generator $[X]$  for each isomorphism class of objects
$X\in \cC$, and we impose the relation
$[X_3] = [X_1]+[ X_2]$, whenever $X_3\cong X_1\oplus X_2$.
When the category $\cC$ possesses  exact sequences, we can define a quotient
$\K(\cC)$ by  imposing the above  relation 
when $X_3$ is an extension of $X_2$ by $X_1$. Although it is not strictly
necessary for our purposes, we will show in the appendix that these
constructions lead to the same groups when applied to $\HS$ and
$\FHS$. Consequently, we will usually drop the label
 ``split'' in the sequel.

 By definition, any additive invariant on $\HS$ or $\FHS$
factors through their\\
Grothendieck groups. 
In particular, this remark applies to
the Poincar\'e polynomial $$P_H(t) = \sum_i \dim H^i\, t^i \in
\Z[t,t^{-1}]$$
and a filtered version of it $$FP_{(H,N)}(t,u) =
\sum_{i,p} \dim (N^p\cap H^i)\, t^iu^p \in \Z[t^{\pm 1}, u^{\pm 1}]$$

We can define two functors between the categories $\HS$ and $\FHS$:
\begin{eqnarray*}
  &\Gamma& : \quad \HS  \longrightarrow  \FHS ~;~  H  \longmapsto
  (H,\F^{\dt})\\
  &\Phi& : \quad \FHS  \longrightarrow  \HS ~;~  (H,N^{\dt}) \longmapsto H
\end{eqnarray*}
where $\F^{\dt}$ is the level filtration on $H$, i.e. $\F^p H$ is the
largest sub Hodge structure of $F^pH$. 
These functors are clearly additive.
 Thus we obtain well-defined group
homomorphisms $\gamma$ and $\phi$, respectively:
\begin{eqnarray*}
  &\gamma& :\quad \K( \HS)  \longrightarrow  \K (\FHS)~;~ \left[ H \right]
  \longmapsto  \left[(H,\F^{\dt})\right]\\
  & \phi&: \quad  \K (\FHS) \longrightarrow  \K (\HS)~;~ \left[(H,N^{\dt})\right]  \longmapsto \left[ H\right]
\end{eqnarray*}

Let $\K (\V_{\C})$ denote the Grothendieck group of the category of
varieties over $\C$ \cite{DenefLoeser1}. A more convenient description
for our purposes is provided by \cite[Theorem 3.1]{Bittner}:
\begin{eqnarray*}
  \K (\V_{\C}) \cong
  \K^{\text{bl}} (\V_{\C}), 
\end{eqnarray*}
where $\K^{\text{bl}} (\V_{\C})$ is the free Abelian group generated
by isomorphism classes $[X]$ of smooth projective varieties subject to
the relation $[\text{Bl}_Z X] - [E] = [X]-[Z]$ for every blow up
$\text{Bl}_Z X$ of $X$ along a smooth closed subvariety $Z \subset X$
with the exceptional divisor $E$.

Let $X$ be a smooth projective variety.  Set
\begin{eqnarray*}
  &&\left[ X \right]_{\HS}:=\sum_i (-1)^i [H^i(X)] \in \K (\HS)\\
  &&\left[ X \right]_{\FHS}:=\sum_i (-1)^i [(H^i(X),N^{\dt})] \in \K (\FHS)
\end{eqnarray*}
where $N^\dt$ is the coniveau filtration.  In the next section, we
will show that these classes depend only on $[X]\in
\K^{\text{bl}}(\V_\C)$.

\section{Coniveau of a blow up}

We use the following notation throughout this section.  Let $X$ be a
smooth projective variety and let $\sigma : \tX=\text{Bl}_Z X \to X$
be the blow up of $X$ along a smooth closed subvariety $Z$ of $X$ of
codimension $\geq 2$. The exceptional divisor $E$ can be identified
with $\mathbb{P}(N_{Z/X})$, where $N_{Z/X}$ is the normal bundle.
Therefore it has a tautological line bundle $\mO_E (1)$.  Let
$r=\codim(Z,X)-1=\dim E - \dim Z$ and let $h= c_1 (\mO_E(1))$.
\begin{eqnarray*}\label{diagram:blowup}
  \begin{aligned}
    \xymatrix{ E ~\ar@{^(->}[rr]^{j} \ar[d]^{\sigma} && \tilde{X}
      \ar[d]^{\sigma}& &
      \tilde{X}-E \ar@{_(->}[ll]_{j_1} \ar[d]^{\sigma}\\
      Z~\ar@{^(->}[rr]^{i} && X && X-Z \ar@{_(->}[ll]_{i_1} }
  \end{aligned}
\end{eqnarray*}
where $i,~i_1,~j,~j_1$ are inclusions.

\begin{lemma}\label{lemma:NpE}
  $$
  N^pH^i(E) = \sigma^*(N^pH^i(Z)) \oplus (h \cup
  \sigma^*(N^{p-1}H^{i-2}(Z))) \oplus \cdots \oplus (h^r \cup\sigma^*(
  N^{p-r} H^{i-2r}(Z))) $$
\end{lemma}

\begin{proof}
  Since $h \in N^1 H^2(E)$, $h^k \in N^k H^{2k}(E)$ for each $ k \geq
  1$. Hence for any $\alpha_k \in \sigma^*(N^{p-k}H^{i-2k}(Z))
  \subseteq N^{p-k}H^{i-2k} (E)$, we have $$
  h^k \cup \alpha_k \in N^k
  H^{2k}(E) \cup N^{p-k} H^{i-2k}(E) \subseteq N^pH^i(E) $$
  by
  \cite[Corollary 1.2]{ArapuraKang} for $1 \leq k \leq r$. Hence this
  and $\sigma^*(N^pH^i(Z)) \subseteq N^pH^i(E)$ give $$
  N^pH^i(E)
  \supseteq \sigma^*(N^pH^i(Z)) \oplus (h \cup
  \sigma^*(N^{p-1}H^{i-2}(Z))) \oplus \cdots \oplus (h^r \cup
  \sigma^*( N^{p-r} H^{i-2r}(Z))) $$

  To show the converse, first note that we have \cite[Proposition
  8.23]{Lewis}
  \begin{eqnarray*}\label{eqn:E}
    H^i (E) = \sigma^*(H^i(Z)) \oplus \left( \bigoplus^r_{k=1} (h^k \cup
      \sigma^*(H^{i-2k}(Z))) \right)  
  \end{eqnarray*}
  Hence for any $\alpha \in N^pH^i(E)$, we can decompose $$
  \alpha =
  \sigma^*(\alpha_0) + (h \cup \sigma^*(\alpha_1)) + (h^2 \cup
  \sigma^*(\alpha_2 )) + \cdots + (h^r \cup \sigma^*(\alpha_r)) $$
  where $\alpha_k \in H^{i-2k}(Z)$ for each $k=0, ...,r$. We will show
  that $\alpha_k \in N^{p-k}H^{i-2k}(Z)$ for each $k=0,...,r$ by using
  descending induction on $k$.

  First note that $$
  \sigma_*(h^i) = \begin{cases} 0 \quad \text{if }
    i < r \\ 1 \quad \text{if } i=r \end{cases} $$
  Therefore, by the
  projection formula and \cite[Theorem 1.1.(1)]{ArapuraKang}, $$
  \alpha_r =\sigma_*(h^r \cup \sigma^*(\alpha_r)) = \sigma_*(\alpha)
  \in N^{p-r} H^{i-2r} (Z) $$
  This shows the claim when $k=r$. Now
  suppose that $\alpha_l \in N^{p-l} H^{i-2l}(Z)$ for $ k +1 \leq l
  \leq r$. We have to show that $\alpha_k \in N^{p-k} H^{i-2k}(Z)$. By
  assumption, for $k+1 \leq l \leq r$ we have $$
  h^l \cup
  \sigma^*(\alpha_l) \in N^l H^{2l}(E) \cup N^{p-l}H^{i-2l}(E)
  \subseteq N^pH^i(E) $$
  Set $$
  \beta_{k+1} = \alpha - \sum^r_{l=k+1}
  (h^l \cup \sigma^*(\alpha_l)) $$
  Then, we have
  \begin{eqnarray*}
    \beta_{k+1}  =
    \sigma^*(\alpha_0) + (h \cup \sigma^*(\alpha_1))+ \cdots + (h^k \cup
    \sigma^*(\alpha_k)) \in N^pH^i(E)
  \end{eqnarray*}
  Then by taking a cup product with $h^{r-k} \in
  N^{r-k}H^{2(r-k)}(E)$, we get
  \begin{eqnarray*}
    h^{r-k} \cup \beta_{k+1} &=& (h^{r-k} \cup \sigma^*(\alpha_0)) +
    (h^{r-k+1} \cup  \sigma^*(\alpha_1)) + \cdots + (h^{r-k+k} \cup
    \sigma^*(\alpha_k)) \\
    &\in& N^{p+(r-k)} H^{i+2(r-k)}(E) 
  \end{eqnarray*}
  Then, $$
  \sigma_*(h^{r-k} \cup \beta_{k+1}) = \sigma_*(h^r \cup
  \sigma^*(\alpha_k)) = \alpha_k \in N^{p-k} H^{i-2k}(Z) $$
  as we
  claimed. Therefore we have $$
  N^pH^i(E) \subseteq
  \sigma^*(N^pH^i(Z)) \oplus (h \cup \sigma^*(N^{p-1}H^{i-2}(Z)))
  \oplus \cdots \oplus (h^r \cup \sigma^*(N^{p-r} H^{i-2r}(Z))) $$
\end{proof}

\begin{cor}\label{cor:NpE}
  Given $\alpha \in N^pH^i(E)$, we can write $$
  \alpha =
  \sigma^*(\alpha_0) + (h \cup \beta) $$
  where $\alpha_0 \in
  N^pH^i(Z)$ and $\beta \in N^{p-1} H^{i-2}(E)$
\end{cor}

The following is well known, but we indicate the proof for lack of a
suitable reference.

\begin{lemma}\label{lemma:ses}
  There is a short exact sequence of pure Hodge structures of weight
  $i$:
  \begin{equation}\label{ses}
    \xymatrix{
      0 \ar[r]& H^i (X) \ar[r]^{\sigma^* + i^* \qquad} &
      H^i(\tilde{X}) \oplus H^i(Z) \ar[r]^{\qquad -j^*+\sigma^*} & H^i(E)
      \ar[r]& 0 
    }
  \end{equation}
\end{lemma}

\begin{proof}

  We have a commutative diagram with exact rows:
  \begin{eqnarray*}
    \xymatrix@C=11pt{
      \cdots \ar[r] & H^i_c (X-Z) \ar[d]^{\cong} \ar[r]&H^i(X) \ar@{^(->}[d]^{\sigma^*} \ar[r]^{i^*} & H^i(Z) \ar[d]^{\sigma^*} \ar[r]^{\partial \qquad }& H^{i+1}_c (X-Z) \ar[d]^{\cong} \ar[r] & H^{i+1}(X) \ar@{^(->}[d]^{\sigma^*} \ar[r]& \cdots\\
      \cdots \ar[r] & H^i_c (\tilde{X}-E) \ar[r]&H^i(\tilde{X}) \ar[r]^{j^*} & H^i(E) \ar[r]^{\tilde{\partial} \qquad }& H^{i+1}_c (\tilde{X}-E) \ar[r] & H^{i+1}(\tilde{X}) \ar[r]& \cdots
    }
  \end{eqnarray*}
  where some of the arrows are injections or isomorphisms as
  indicated. The lemma now follows from a straight forward diagram
  chase.
\end{proof}

\begin{lemma}\label{lemma:sesNp} 
  The sequence of the previous lemma is exact in $\FHS$, i.e.  $$
  \xymatrix{ 0 \ar[r]& N^pH^i (X) \ar[r] & N^pH^i(\tilde{X}) \oplus
    N^pH^i(Z) \ar[r] & N^pH^i(E)\ar[r]&0 } $$
  is an exact sequence for
  all $p$.
\end{lemma}

\begin{proof}
  Consider the short exact sequence (\ref{ses}) and note that
  \begin{eqnarray*}
    (\sigma^* +i^*)(N^pH^i(X))\subseteq N^pH^i(\tX) \oplus N^pH^i(Z),\\
    (-j^*+\sigma^*)(N^pH^i(\tX) \oplus N^pH^i(Z))\subseteq N^pH^i(E)
  \end{eqnarray*}
  since all maps preserve the coniveau.

  We will check exactness in the middle. It suffices to show that $$
  \im(\sigma^* +i^*)|_{N^pH^i(X)} \supseteq \ker (-j^* +
  \sigma^*)|_{N^pH^i(\tX) \oplus N^pH^i(Z)}, $$
  since the reverse
  inclusion follows from (\ref{ses}).  Let $(\beta, \gamma) \in
  N^pH^i(\tX) \oplus N^pH^i(Z)$ such that $(\beta, \gamma) \in
  \ker(-j^* +\sigma^*)$. Then, from the exact sequence (\ref{ses}),
  there is $\alpha \in H^i(X)$ such that $(\sigma^*
  +i^*)(\alpha)=(\beta, \gamma)$. In particular, we have $$
  \sigma^*(\alpha) =\beta \in N^pH^i(\tX) $$
  Since $\sigma: \tX \to X$
  is a birational map, we have $$
  \alpha = (\sigma_* \circ
  \sigma^*)(\alpha) =\sigma_*(\beta) \in \sigma_*(N^pH^i(\tX))
  \subseteq N^pH^i(X) $$
  Hence we have $(\beta, \gamma)=(\sigma^*+
  i^*)(\alpha)\in (\sigma^*+i^*)(N^pH^i(X))$.

  We will check surjectivity of $ (-j^* + \sigma^*)|_{N^pH^i(\tX)
    \oplus N^pH^i(Z)}$. Let $\alpha \in N^pH^i(E)$.  Then by
  Corollary~\ref{cor:NpE} we can decompose $$
  \alpha =
  \sigma^*(\alpha_0) + (h \cup \beta) $$
  where $\alpha_0 \in
  N^pH^i(Z)$ and $\beta \in N^{p-1} H^{i-2} (E)$.  Note the
  composition $$
  \xymatrix{ H^{i-2} (E) \ar[r]^{j_*} & H^i(\tX)
    \ar[r]^{j^*} & H^i(E) } $$
  is same as cupping with $[E]|_E$. Since
  $\mO(E)|_E= \mO_{E}(-1)$ and $h = c_1 (\mO_E(1))$, we have $$
  (j^*
  \circ j_*)(\beta)= [E]|_{E} \cup \beta =-h \cup \beta $$
  Thus, we
  have $$
  \alpha = \sigma^*(\alpha_0) - j^*(j_*(\beta))$$
  and
  $j_*(\beta) \in j_*(N^{p-1}H^{i-2}(E))\subseteq N^p H^i(\tX)$ and
  this shows that the map $(-j^* +\sigma^*)|_{N^pH^i(\tX) \oplus
    N^pH^i(Z)}$ is surjective.
\end{proof}

Let $f: X \to Y$ be a morphism of smooth projective varieties. Then
$f^*$ preserves the coniveau and induces maps $\bar{f^*}:Gr^p_N H^i(Y)
\to Gr^p_N H^i(X)$.

\begin{cor}\label{cor:sesGr} The following sequence is exact:
  $$
  \xymatrix{ 0 \ar[r]& Gr^p_N H^i(X) \ar[r] &Gr^p_N H^i(\tilde{X})
    \oplus Gr^p_N H^i(Z) \ar[r] &Gr^p_N H^i(E) \ar[r]& 0 } $$
\end{cor}

\begin{proof}
  This follows from lemma \ref{lemma:exactness}.
\end{proof}

We want to construct well-defined morphisms from
$\K^{\text{bl}}(\V_{\C})$ to $\K(\HS)$ and from
$\K^{\text{bl}}(\V_{\C})$ to $\K(\FHS)$. First we need the following
lemmas.

\begin{lemma}\label{lemma:HS}
  The relation $$
  [\text{Bl}_Z X]_{\HS} - [E]_{\HS} = [X]_{\HS}-
  [Z]_{\HS} $$
  holds. Hence, we have a well-defined group homomorphism
  \begin{eqnarray*}
    \K^{\text{bl}}(\V_{\C}) \longrightarrow \K(\HS)~;~ \left[ X \right]\mapsto  [X]_{\HS}
  \end{eqnarray*}

\end{lemma}

\begin{proof}
  By Lemma~\ref{lemma:ses} we have a short exact sequence of pure
  Hodge structures:
  \begin{eqnarray*}
    \xymatrix{
      \qquad 0 \ar[r]& H^i (X) \ar[r] &
      H^i(\tilde{X}) \oplus H^i(Z) \ar[r]& H^i(E)
      \ar[r]& 0 
    }
  \end{eqnarray*}
  Then by the definition of $\K(\HS)$, we have $$
  [H^i(X)] +
  [H^i(E)]=[H^i(\tilde{X}) \oplus H^i(Z)] =[H^i(\tilde{X})]+[H^i(Z)]
  $$
  and the result follows immediately by taking alternating sums.
\end{proof}

We define $$\lambda: \K^{\text{bl}}(\V_{\C}) \longrightarrow
\K(\FHS)$$
as the composition of the map constructed in the above
lemma with $\gamma$.

\begin{lemma}\label{lemma:FHS}
  The relation $$
  [\text{Bl}_Z X]_{\FHS} - [E]_{\FHS} = [X]_{\FHS}-
  [Z]_{\FHS} $$
  holds. Hence we have a well-defined group homomorphism
  \begin{eqnarray*}
    \nu: \K^{\text{bl}}(\V_{\C}) \longrightarrow \K(\FHS)~;~\nu(\left[ X \right])= [X]_{\FHS}
  \end{eqnarray*}

\end{lemma}

\begin{proof}
  We showed in Lemma~\ref{lemma:sesNp} that $$
  0 \to [(H^i(X),
  N^{\dt})] \to [(H^i(\tilde{X}) \oplus H^i(Z), N^{\dt})]\to [(H^i(E),
  N^{\dt})] \to 0 $$
  is exact in $\FHS$ where
  $N^p(H^i(\tilde{X})\oplus H^i(Z))=N^pH^i(\tilde{X}) \oplus
  N^pH^i(Z)$. The rest of the argument is exactly the same as above.
\end{proof}

\begin{rmk}
The proof of  Lemma~\ref{lemma:sesNp} shows that the above  
sequence is split exact. So in particular,  as the referee has
pointed to us, we can construct the homomorphism
$$ \nu: \K^{\text{bl}}(\V_{\C}) \longrightarrow \K^{split}(\FHS)$$
directly, without appealing to the results of the appendix.
\end{rmk}

\section{Main theorem}

\begin{thm}\label{thm:GHC}
  Let $X$ be a smooth projective variety. Then, the following
  statements are equivalent:
  \begin{enumerate}
  \item[1.]  GHC holds for $X$.
  \item[2.] $[X]\in \ker(\nu-\lambda)$.
  \item[3.] The equality $FP_{[X]_\FHS}(t,u)=
    FP_{\gamma([X]_\HS)}(t,u)$ of filtered Poincar\'e polynomials
    holds.
  \end{enumerate}
\end{thm}

\begin{proof}
  It is clear that GHC for $X$ implies $[X]_{\FHS} = \gamma
  ([X]_{\HS})$, or equivalently that $ [X]\in \ker(\nu-\lambda)$.

  Suppose $[X]\in \ker(\nu-\lambda)$.  Then $$
  [X]_{\FHS} = \gamma
  ([X]_{\HS}) $$
  i.e.
  \begin{eqnarray*}
    \sum_i (-1)^i [(H^i(X),N^{\dt})] &=& \sum_i (-1)^i [(H^i(X),\F^{\dt})]
  \end{eqnarray*}
  Taking the filtered Poincar\'e polynomial $FP(t,u)$ of both sides
  yields the third statement.

  Assume the equality in 3.  The coefficient of $t^i$ on the left is
  $$(-1)^i\sum_p \dim N^pH^i(X)u^p,$$
  while on the right it is
  $$(-1)^i\sum_p \dim \F^pH^i(X)u^p$$
  The equality of these
  expressions forces $ N^pH^i(X)=\F^pH^i(X)$ in this case.
\end{proof}

\begin{rmk}\label{remark:GHC}
  The coefficient of $t^iu^p$ in $FP_{[X]_\FHS}(t,u)$ is $(-1)^i\dim
  N^pH^i(X)$, and the coefficient of $t^iu^p$ in
  $FP_{\gamma([X]_\HS)}(t,u)$ is $(-1)^i\dim \F^p H^i(X)$.  Therefore
  $\text{GHC}(H^i(X),p)$ holds precisely when these coefficients
  coincide.
\end{rmk}

Let $\LL'=[\A^1_{\C}] \in \K(\V_{\C})$ denote the Lefschetz object.
Under the isomorphism $$\K(\V_{\C}) \cong \K^{\text{bl}}(\V_{\C})$$
$\LL'$ maps to $\LL=[\PP^1] - [\infty]$. Let $$\KK= [(\Q(-1),
N^{\dt})]$$
where $N^\dt$ is the filtration determined by $N^1/N^2 =
\Q(-1)$.  Then $$
\lambda(\LL) =\nu(\LL)= \KK$$
We have a product on
the category $\FHS$ which is just the tensor product filtered by $$
N^p (H_1\otimes H_2) = \bigoplus_{r+s=p} N^r H_1 \otimes N^s H_2 $$
This gives a commutative ring structure on $\K(\FHS)$.  The group
$\K^{\text{bl}}(\V_{\C})$ also has a commutative ring structure
induced by the product of varieties.  It is not clear whether
$\lambda$ or $\nu$ are ring homomorphisms, however we do have:

\begin{lemma} 
  For any $\eta \in \K^{\text{bl}}(\V_{\C})$, $$\lambda(\LL\cdot \eta)
  = \KK\cdot \lambda(\eta)$$
  $$\nu(\LL\cdot \eta) = \KK\cdot \nu(
  \eta)$$
\end{lemma}

\begin{proof}
  Under the K\"unneth isomorphism, we have $$N^pH^i(\PP^1\times X)
  \cong N^{p-1}H^{i-2}(X)(-1)\oplus N^pH^i(\{\infty\}\times X),$$
  and
  a similar statement holds for $\F^p$. The lemma is an immediate
  consequence.
\end{proof}

The element $\KK$ is invertible, and the above identities guarantee
that $\lambda$ or $\nu$ factor through the localization
$\M=\K^{\text{bl}}(\V_{\C})[\LL^{-1}]$.

Recall (\cite{DenefLoeser1}) that there is a decreasing filtration
$F^{\dt}\M$ on $\M$, where $F^m \M $ is the subgroup of $\M$ generated
by $\{[X] \cdot \LL^{-j}~|~ \dim X-j \leq -m \}$.  Let $\hM$ be the
completion of the ring $\M$ with respect to the filtration
$F^{\dt}\M$.  A similar filtration (compare \cite{Loo}) can be defined
on $\K(\FHS)$ by replacing dimension by weights. More precisely, let
$L^m\K(\FHS)$ be the subgroup generated by $\{[(H,N^\dt)] \mid \max
\{i\mid H^i\not=0\}\leq -m\}$.  We denote the completion of $\K(\FHS)$
with respect to $L^{\dt}$ by $\hN$.  The weights of the Hodge
structure on cohomology of a smooth projective variety of dimension
$\le d$ are bounded by $2d$. Therefore the induced filtrations
$\lambda(F^{\dt})$ and $\nu(F^{\dt})$ are cofinal with a subfiltration
of $L^{\dt}$.  It follows that we have a commutative diagram
\begin{eqnarray*}\label{diagram:inj}
  \begin{aligned}
    \xymatrix{
      \K^{\text{bl}} (\V_{\C}) \ar[rr]^{\qquad f} \ar[d]^{\alpha} && \hM \ar[d] \\
      \K(\FHS)\ar[rr]^{\qquad \tau} && \hN }
  \end{aligned}
\end{eqnarray*}
where $\alpha$ denotes either $\lambda$ or $\nu$, and $f,\tau$ are the
canonical ones.

\begin{lemma}\label{lemma:inj}
  The homomorphism $\K(\FHS)\to \Z[t^{\pm 1}, u^{\pm 1}]$ given by the
  filtered Poincar\'e polynomial factors through the image of $\tau$.
\end{lemma}

\begin{proof}
  Let $\eta \in \bigcap_m L^m $.  Then for each $m\in \Z$, $\eta$ can
  be expressed as a linear combination of classes of filtered Hodge
  structures of weight at most $-m$. Thus the degree of $FP_\eta(t,u)$
  in $t$ is bounded above by $-m$, for all $m\in \Z$. This is
  impossible unless $FP_\eta(t,u)=0$.
\end{proof}

\begin{cor}[to theorem]\label{cor:GHC}
  Let $X,Y$ be smooth projective varieties such that $X,Y$ define the
  same classes in $\hM$. Then GHC holds for $X$ if and only if GHC
  holds for $Y$.
\end{cor}

\begin{proof}
  If the images of $[X]$ and $[Y]$ coincide in $\hM$ then they
  coincide in $\hN$. Therefore their filtered Poincar\'e polynomials
  coincide.  The conclusion is now an immediate consequence of
  Theorem~\ref{thm:GHC}.
\end{proof}

From Remark~\ref{remark:GHC}, we obtain:

\begin{cor}
  Let $X,Y$ be smooth projective varieties such that $X,Y$ define the
  same classes in $\hM$. Then for each $i$ and $p$,
  $\text{GHC}(H^i(X),p)$ holds if and only if $\text{GHC}(H^i(Y),p)$
  holds.
\end{cor}

The proof of the next corollary depends on the motivic integration
theory of Kontsevich, Denef and Loeser. See \cite{DenefLoeser1},
\cite{DenefLoeser2}, and \cite{Loo} for an introduction to these
ideas.

\begin{cor}\label{cor:Keq} 
  Let $X,Y$ be $K$-equivalent smooth projective varieties, i.e. there
  is a smooth projective variety $Z$ and birational maps $\pi_1:Z \to
  X$ and $\pi_2:Z \to Y$ such that $\pi_1^* K_X = \pi^*_2 K_Y$, where
  $K_X$ (respectively $K_Y$) is the canonical divisor on
  $X$ (respectively $Y$).  Then $\text{GHC}(H^i(X),p)$ holds if and
  only if $\text{GHC}(H^i(Y),p)$ holds.
\end{cor}

\begin{proof}
  It is enough to show that $X$ and $Y$ define the same class in
  $\hM$. This follows from the $K$-equivalence assumption by a
  standard application of motivic integration theory, see
  \cite{loeser} or \cite{veys}. For convenience of the reader, we
  reproduce the argument.  By the change of variables formula
  \cite[Lemma 3.3]{DenefLoeser1}, we have
  \begin{eqnarray*}
    f([X]) = \int_{\cL(X)} d \mu_X &=& 
    \int_{\cL(Z)} \mathbb{L}^{-{ord_t\pi_1^*\omega_{X}}} d \mu_Z \\
    &=& \int_{\cL(Z)} \mathbb{L}^{-{ord_t\pi_2^*\omega_{Y}}} d \mu_Z =\int_{\cL(Y)} d \mu_Y =f([Y])
  \end{eqnarray*}
  where $f: \K^{\text{bl}}(\V_{\C}) \to \hM$ is the canonical map,
  $\cL(X),\cL(Y)$ are the arc spaces, $\omega_X, \omega_Y $ are the
  canonical sheaves, and $d\mu_X,d\mu_Y$ are the motivic measures.
  Hence the corollary follows from Corollary~\ref{cor:GHC}.
\end{proof}

\begin{cor}\label{cor:CalabiYau}
  Let $X,Y$ be birational Calabi-Yau varieties. Then
  $\text{GHC}(H^i(X),p)$ holds if and only if $\text{GHC}(H^i(Y),p)$
  holds.
\end{cor}

\begin{proof}
  Since $K_X=K_Y=0$, it follows from Corollary~\ref{cor:Keq}.
\end{proof}

\begin{appendix}
\section{Grothendieck groups of filtered categories}

We recall that an exact category consists of an
 additive category $\cC$, together with a distinguished class
 of diagrams $$0\to A\to B\to C\to 0$$
 called exact sequences  satisfying appropriate conditions \cite{Quillen}.
For example, any additive category $\cC$ can be made 
exact by taking the class of exact sequences to be isomorphic to the 
class of split sequences 
$$0\to A\to A\oplus C\to C\to 0.$$
An Abelian category gives another example of an exact category, where
exact sequences have the usual meaning.
If the category is also semisimple, then this exact structure
coincides with the split structure above. This remark applies to $\HS$. 

 Let $\cC$ be an Abelian category, and let
$\FcC$ (respectively $\GcC$) denote the category of filtered
(respectively graded) objects in $\cC$. The category $\GcC$
is Abelian, but  $\FcC$ is generally not. However $\FcC$  has a natural exact
structure \cite[1.1.4]{bbd} given as follows.
 We have a functor $$(H,N^{\dt}) \longmapsto \bigoplus_p N^p H$$
 from $\FcC$ to $\GcC$.  
We declare a sequence $$
 0 \to (H_1, N^{\dt}) \to (H_2, N^{\dt}) \to (H_3,
 N^{\dt}) \to 0 $$
 in $\FcC$ to be {\em exact} if and only if its
 image in $\GcC$ is exact. For the record, we note the following
 alternative formulation, which is perhaps more common:

\begin{lemma}\label{lemma:exactness}
  The sequence $$
  0 \to (H_1, N^{\dt}) \to (H_2, N^{\dt}) \to (H_3,
  N^{\dt}) \to 0 $$
  is exact if and only if the following sequence is
  exact in $\GcC$: $$
  0 \to \bigoplus_p N^pH_1 /N^{p+1}H_1
  \longrightarrow \bigoplus_p N^pH_2 /N^{p+1}H_2 \longrightarrow
  \bigoplus_p N^pH_3 /N^{p+1}H_3 \to 0 $$
\end{lemma}

\begin{proof}
  This is a straight forward application of the Snake lemma and
  induction.
\end{proof}

\begin{lemma}\label{lemma:FHSisExact}
  The category $\FcC$ with the above notion of exact sequence is an
  exact category.
\end{lemma}

Given an exact category $\cC$, its Grothendieck
group $\K(\cC)$ is given by generators $[M]$, with $ M\in \cC$,
and  relations $[M_2]=[M_1]+[M_3]$ for every exact sequence $0\to
  M_1\to M_2\to M_3\to 0$.
Let us denote the  Grothendieck group for $\cC$ with its split
exact structure  by $\K^{split}(\cC)$.
We see immediately that, $\K^{split}(\cC) \cong \K(\cC)$ if $\cC$ is
Abelian and semisimple. This is, in particular, the case for $\HS$.
We have a homomorphism $\K^{split}(\FHS)\to \K(\FHS)$, which
 is also an isomorphism by:

\begin{lemma}\label{lemma:Fsemisimple}
If $\cC$ is semisimple Abelian, then
any exact sequence in $\FcC$ is split exact.
\end{lemma}

\begin{proof}
It is enough to check that  given an exact sequence in $\FcC$
\begin{equation}\label{seq:exact}
0 \to (H_1, N^{\dt}) \to (H_2 , N^{\dt}) \stackrel{f}{\to} (H_3, N^{\dt}) \to 0,
\end{equation}
there is a splitting $s:(H_3, N^{\dt}) \to (H_2, N^{\dt})$ for $f$. 
First note that by semisimplicity of $\cC$, we have a noncanonical decomposition 
$$
N^pH_i \cong N^{p+1}H_i \oplus Gr^p_N H_i
$$ 
where $Gr^p_N H_i = N^pH_i/N^{p+1}H_i$, $i=2,3$. We define $s_p= s|_{N^{p}H_3}$ by descending induction on $p$ :
Note that the exact sequence (\ref{seq:exact}) induces exact sequences in $\cC$
$$
0 \longrightarrow N^{p+1}H_1 \longrightarrow N^{p+1}H_2  \stackrel{f_{p+1}}{\longrightarrow} N^{p+1}H_3 \longrightarrow 0 
$$
$$
0 \longrightarrow Gr^p_N H_1 \longrightarrow Gr^p_N H_2 \stackrel{\bar{f}_p}{\longrightarrow} Gr^p_N H_3 \longrightarrow 0
$$
where $f_{p+1} = f|_{N^{p+1}H_2}$ and $\bar{f}_p$ is the induced
map. By induction and semisimplicity 
of $\cC$, there are splittings 
$s_{p+1}: (N^{p+1}H_3, N^\dt\cap N^{p+1}H_3)  \to (N^{p+1}H_2
N^\dt\cap N^{p+1}H_2)$ 
and $t_p : Gr^p_N H_3 \to Gr^p_N H_2 $ for $f_{p+1}$ and $\bar{f}_p$, respectively. Set 
$$
s_p = s_{p+1} + t_p : N^pH_3 \longrightarrow N^p H_2
$$
Then $s_p$ gives a well-defined splitting for $f_p$ and hence we have
a splitting $s=s_0$ for $f$. This completes the proof of the Lemma. 
\end{proof}

\begin{cor}\label{cor:K0split}
  $\K^{split} (\FHS) \cong \K (\FHS)$
\end{cor}

\end{appendix}

 \end{document}